\input amstex
\documentstyle{amsppt}

\document

\topmatter
\title
Holomorphic horospherical duality "sphere-cone"
\endtitle

\author Simon Gindikin \endauthor
\dedicatory To  Gerrit van Dijk with warmest feelings
\enddedicatory

\address Departm. of Math., Hill Center, Rutgers University,
110 Frelinghysen Road, Piscataway, NJ 08854-8019
\endaddress
\email gindikin\@math.rutgers.edu \endemail \abstract We describe
a construction of complex geometrical analysis which corresponds
to the classical theory of spherical harmonics\endabstract
\endtopmatter
I believe that the connection of harmonic analysis and complex
analysis has an universal character and is not restricted by the
case of complex homogeneous manifolds. It looks as a surprise that
such a connection exists and though it is quite natural for finite
dimensional representations and compact Lie groups \cite
{Gi00,Gi02}. In this note we describe the complex picture which
corresponds to harmonic analysis on the real sphere. The basic
construction is a version of horospherical transform which in this
case is a holomorphic integral transform between holomorphic
functions on the complex sphere and the complex spherical cone.
This situation looks quite unusual from the point of view of
complex analysis and I believe presents a serious interest also in
this setting. It can be considered as a version of the Penrose
transform, but in a purely holomorphic situation when there is
neither cohomology nor complex cycles.

\subhead Geometrical picture\endsubhead Let
$$\Delta(z)=z_1^2+\dots +z_{n+1}^2$$ be the quadratic form in
$\Bbb C^{n+1}_z$. Let $\Bbb C^{n+1}_z$ and $\Bbb C^{n+1}_\zeta$
are dual spaces relative to the form $$\zeta\cdot
z=\zeta_1z_1+\cdots +\zeta_{n+1}z_{n+1}.$$ Let the complex sphere
$\Bbb CS$=$\Bbb CS^n \subset \Bbb C^{n+1}_z$ be defined by the
equation $$\Delta(z)=1$$ and the cone $\Bbb C\Xi$=$\Bbb C\Xi^n
\subset \Bbb C^{n+1}_\zeta$ be defined as
$$\Delta(\zeta)=0, \zeta\neq0.$$ Both these complex manifolds are
homogeneous relative to $SO(n+1,\Bbb C)$, but we want to describe
our objects in the language of complex geometry without an appeal
to groups.

We consider on $\Bbb CS$ the family of sections $E(\zeta),\zeta
\in \Bbb C\Xi,$ by hyperplanes $$\zeta\cdot z=1.$$ Let us call
$E(\zeta)$ horospheres. They are paraboloids. Correspondingly on
$\Bbb C\Xi$ we consider the hyperplane sections $L(z),z\in \Bbb
CS,$ by the hyperplanes with the same equations. They are
hyperboloids. These dual families of complex submanifolds are the
basic element of our geometrical picture.

There is one more essential element of the picture: we consider
the family of real spheres $S(u), u\in U,$ which are real forms of
the complex sphere $\Bbb CS$. The parametrical space $U$ is the
homogeneous space $SO(n+1,\Bbb C)/SO(n+1)$, but for us it is
essential that $S(u)$ are totally real cycles of the (maximal)
dimension $n$. They are all mutually homotopic, but they are not
homological to zero. \subhead Horospherical Cauchy
transform\endsubhead The complex sphere $\Bbb CS$ is the Stein
manifold; the complex cone $\Bbb C\Xi$ is the holomorphicaly
separated complex manifold. Thus there are many holomorphic
functions on both of them. We will construct some remarkable
integrable operators between the spaces of holomorphic functions
on them.

We will denote through $[a_1, a_2,\dots, a_{n+1}]$ the determinant
of the matrix with the columns $a_1, a_2,\dots, a_{n+1}$ some of
which can be 1-forms. We expand such determinants from left to
right using the exterior product for the multiplication of
1-forms. We will write $a^{\{k\}}$ if a column $a$ repeats $k$
times. Let
$$\omega =[z,dz^{\{n\}}]=n!\sum_{1\leq j\leq (n+1)}(-1)^{j-1}z_j\bigwedge_{i\neq j}dz_i.$$  It is an
invariant holomorphic $n$-form on $\Bbb CS.$

For a holomorphic function $f(z)\in \Cal O(\Bbb CS)$ let$$ \hat
f(\zeta) = \int_ S f(z) \frac \omega {1 -\zeta\cdot z}.$$ Here $S$
is any real sphere $S(u)$ which does not intersect the horosphere
$E(\zeta)$. Such cycles exist and they are homological. The
integral is independent of the choice $S(u)$ since the integrand
is closed as a holomorphic form of maximal degree. Let us remark
that this form has the singularity on the horosphere $E(\zeta)$
and that we integrate it along the compact manifold. The result is
a holomorphic function on $\Bbb C\Xi$. This operator $f\rightarrow
\hat f$ acting from $\Cal O(\Bbb CS)$ to $\Cal O(\Bbb C\Xi)$ we
will call {\it the horospherical Cauchy transform}.

For a real sphere $S(u)$ let $\Xi(u)$ be the set of such $\zeta
\in \Xi$ that $S(u)\bigcap E(\zeta) =\emptyset$. It is simple to
check \cite {Gi04} that, if $S(0)=S$ is the sphere $\{y=\Im
z=0\}$, then $\Xi(0)=\Xi_S$ is the domain described by the
conditions
$$\Delta(\xi)=\Delta(\eta)<1.$$ These domains $\Xi(u)$ for
different $u$ can be obtained  by the action of elements of the
group $G_\Bbb C=SO(n+1,\Bbb C)$. Dual connected sets $U(\zeta)$
are described similarly.

\subhead Dual horospherical Cauchy transform \endsubhead Let us
consider on $L(z)$ the normalized holomorphic forms of the maximal
degree $n-1$: $$\nu_z(d\zeta)=
\frac{[\lambda,\zeta,d\zeta^{\{n-1\}}]}{\lambda \cdot z},$$ where
$\lambda$ is any vector such that $\lambda \cdot z\neq 0$. The
restriction of the form $\nu_z$ on $L_z$ is independent of a
choice of $\lambda$: up to a constant factor it is the residue of
the form $$\frac {\omega (d\zeta)}{1-\zeta \cdot z}$$ on $L_z$. We
can take, in particular, $\lambda = z$ and then $$ \nu
_z(d\zeta)=[z,\zeta,d\zeta^{\{n-1\}}].$$ Let us put
$$ \check F(z)=\int _{L_\Bbb R(z)} F(\zeta)\nu_z(d\zeta), F\in
\Cal O(\Bbb C\Xi), $$  Here $L_\Bbb R(z)$ is any cycle in $L(z)$
which is its real form (similar to $S(u)$). We have an operator
from $\Cal O(\Bbb C\Xi)$ to $\Cal O(\Bbb C S)$ which we will call
{\it the dual horospherical Cauchy transform}.

Let us compare the two transforms which we defined. They are very
similar indeed. To see this we remark that in the definition of
the dual transform we can replace the integral by the integral of
the $n$-form $$\frac {F(\zeta) [\zeta, d\zeta^{\{n\}}]}{1-z\cdot
\zeta}$$ along a $n$-dimensional cycle contractible to ${L_\Bbb
R(z)}$. The integral in the definition is the residue of this
integral. This version shows that $\check F$ is holomorphic on $z$
since we can not change the cycle for small changes of $z$. So the
definitions are similar but there is an essential difference
originating in a difference in the geometry of $\Bbb C S$ and
$\Bbb C \Xi$: we can take the residue of $\omega$ on the
horosphere $E(\zeta)$ but we can not contract the cycle $S(u)$ on
$E(\zeta)$ and on the horosphere there are no appropriate cycles.
In a sense $\hat f$ characterizes $f$ on "infinity".

\subhead The inversion formula \endsubhead We want to find a
Radon's type inversion formula for the horospherical transform.
This formula will combine the dual transform and a remarkable
differential operator $\Cal L$ on $\Bbb C\Xi$. Let $$ D=\zeta
\cdot \frac {\partial}{\partial \zeta}= \zeta_1
\frac{\partial}{\partial \zeta_1}+\cdots +\zeta_{n+1}
\frac{\partial}{\partial \zeta_{n+1}}$$ be the operator along
generators of the cone $\Bbb C\Xi$. The operator $\Cal L$ will be
a polynomial $\Cal L(D)$ of $D$. We extend functions $F(\zeta)$ as
homogeneous functions $F(\zeta,p),p\in \Bbb C,$ of degree -1. It
corresponds to the homogeneous coordinates of the sections
($\zeta\cdot z=p$). Then
$$\Cal L=c(\frac {n-1} 2 \frac {\partial ^{(n-2)}}{\partial p^{(n-2)}} -
2\frac {\partial ^{(n-1)}}{\partial p^{(n-1)}})|_{p=1}, \quad
c=\frac n {(-2\pi i)^n}.$$ This operator is the polynomial $\Cal
L(D)$ of $D$.

\proclaim {Theorem} There is an inversion formula $$f=(\Cal L \hat
f)^{\vee},\quad f\in \Cal O(\Bbb CS).$$ \endproclaim It is
sufficient (using the closeness of forms) to check the formula in
one point $z\in \Bbb CS$ and for a specific choice of the cycles
$S(u)$ and $L_\Bbb R(z)$. Let $z=x$ be a point of the real sphere
$S=S(0)$ and $$L_\Bbb R(x)=\{ \zeta= x+i\eta, x\cdot
\eta=\eta\cdot\eta=0\}.$$ The inversion formula in this
specification was proved in \cite {Gi04} for arbitrary continuous
functions on the real sphere $S$.

We can consider this inversion formula as the integral formula of
Cauchy-Fantappie type on $\Bbb CS$ (cf.\cite {GH90}).

\subhead The connection with spherical harmonics \endsubhead All
our constructions so far can be described only in the language of
the geometry of two families of complex submanifolds. They all are
invariant relative to $SO(n+1;\Bbb C)$ but we can avoid the group
language. Now in our considerations will appear a group but this
group will be Abelian: the Abelian group $\Bbb C^\times$ acts on
$\Bbb C\Xi$ by the multiplications $\zeta\mapsto c\zeta, c\in \Bbb
C\setminus 0$. This action commutates with the action of
$SO(n+1;\Bbb C)$. Let us decompose the space of holomorphic
functions $\Cal O(\Bbb C\Xi)$ on invariant subspaces $$\Cal O(\Bbb
C\Xi)=\bigoplus _{n\geq 0} \Cal O_n(\Bbb C\Xi).$$ Functions in
$\Cal O_n(\Bbb C \Xi)$ are homogeneous polynomials of degree $n$
and they can be interpreted as sections of line bundles on the
projectivization $F$ of the cone $\Bbb C\Xi$ (the flag manifold).
In $\Cal O_n(\Bbb C\Xi)$ there are realized the irreducible
(finite dimensional) representations of $SO(n+1;\Bbb C)$.

Let us call {\it the holomorphic Fourier transform} $\tilde
f(n;\zeta)$ the composition of the horospherical Cauchy transform
and the projection on $\Cal O_n$.  By the direct decomposition of
the kernel $1/(1-\zeta\cdot z)$ in the geometrical series we
obtain
$$\tilde f(n;\zeta)=\int _{S(u)} f(z)(\zeta\cdot z)^n\omega.$$

The differential operator $\Cal L(D)$ on $\Cal O_n$ is the
multiplication on the $\Cal L(n)$. As a result we have the
inversion formula for the spherical Fourier transform $$
f(z)=\sum_ {n\geq 0} \Cal L(n) f_n(z),\quad f_n(z)=\int_{L_\Bbb
R(z)}\tilde f(n;\zeta)\nu_z(d\zeta).$$ The preimages of $\Cal
O_n(\Bbb C\Xi)$ relative to the horospherical transform are
subspaces of spherical harmonics $\Cal O_n(\Bbb CS)$; $f\mapsto
f_n$ are the projectors on them. The integral operator $f(n;\zeta)
\mapsto f_n(z)$ can be considered as an analogue of the Poisson
integral. The quadric $F$ in a sense plays the role of "the
complex boundary" of $\Bbb CS$. The connection of spherical
harmonics and homogeneous polynomials on the cone $\Bbb C\Xi$ goes
back to Maxwell.

\subhead Holomorphic extensions of horospherical transform
\endsubhead We can remove in the definition of the horospherical
transform $\hat f(\zeta) $ the condition $\Delta (\zeta)=0$. We
need only the existence of a cycle $S(u)$ which does not intersect
the section by $\zeta\cdot z=1$. Correspondingly, we can extend
the definition of the dual horospherical transform for such $z$
that the hyperplane $ z\cdot \zeta=1$ intersects the cone $\Bbb
C\Xi$ on a hyperboloid. The extended holomorphic functions will
satisfy to some differential equation. The direct differentiation
of the integral representation of $\check F(z)$ shows that it is
harmonic: $$\Delta( \frac \partial {\partial z})\check F(z)=0.$$
The combination with the inversion formula gives the harmonic
extension of holomorphic functions $f(z)$ from $\Bbb CS$. The
elements of $\Cal O_n(\Bbb CS)$ extend as harmonic polynomials.

To describe the extension of $\hat f(\zeta)$ it is convenient to
extend them as homogeneous functions $\check f(\zeta, p)$. Then we
have $$ [\Delta( \frac \partial {\partial \zeta})- \frac {\partial
^2}{\partial p^2}] \hat f(\zeta,p)=0.$$ Of course using the
homogeneity of $\hat f$ we can eliminate $p$.

\subhead Unitary restriction\endsubhead The principle of the
unitary restriction means that irreducible finite dimensional
representations of complex semisimple Lie groups give the complete
system of such representations for compact forms of these groups.
In our geometrical picture it corresponds to the possibility to
restrict the construction of the horospherical Cauchy transform on
the real sphere $S$ (indeed we already used this possibility in
the opposite direction). More precisely, let us consider for the
sphere $S\subset \Bbb CS$ the dual domain $\Xi_S\subset \Bbb C\Xi$
of such $\zeta$ that the horosphere $E(\zeta)$ does not intersect
$S$. As we remarked (with the reference \cite {Gi04}) this domain
is described by the condition
$$ \Delta(\xi)=\Delta (\eta)<1.$$ Correspondingly, we can construct
the duality between the real spheres $S(u)$ and the domains
$\Xi(u)$ in $\Bbb C\Xi$. I believe that it is remarkable that the
natural dual object to the real sphere is a complex domain.

For functions $f\in C(S)$ we can define using our definition the
horospherical Cauchy transform $\hat f(\zeta), \zeta \in \Xi_S.$
We can take the boundary values of $\hat f$ and apply the
inversion formula. The crucial circumstance is that submanifold
$L(x),x\in S,$ does not intersect the domain $\Xi_S$ but intersect
the boundary on the cycle $$L_\Bbb R(x)= \{ \zeta=x+i\eta, x\cdot
\eta=\eta\cdot\eta=0\}.$$ In such a way the boundary $\partial
\Bbb (\Bbb C\Xi)$ fibers over $S$ with the fibers $L_\Bbb R(x)$.
The inversion formula holds in this situation \cite {Gi04}.

Since the domain $\Xi_S$ is invariant relative to the
multiplication on the circle ($\zeta \mapsto \exp(i\theta)\zeta
$), the decomposition in the Fourier series survives in this
domain. Functions in invariant subspaces are automatically
holomorphic in $\Bbb C\Xi$ and $\Cal O_n(\Xi_S)=\Cal O_n(\Bbb
C\Xi)$. Functions in dual subspaces on $S$ also will be extended
holomorphically on $\Bbb CS$. It reflects the elliptic nature of
the problem (these functions are eigenfunctions of invariant
differential operators on $S$ which are elliptic). The holomorphic
Fourier transform and its inversion are preserving for $f\in
C(S)$.

It is natural to discuss functional spaces on $S$ for which these
constructions make sense. We can use our definition of the
horospherical Cauchy transform for distributions on $S$. Moreover,
we can consider the horospherical transform of hyperfunctions. So
we consider the space of hyperfunctions
$$\operatorname {Hyp}(S) =H^{(n-1)}(\Bbb CS\backslash S, \Cal
O), $$ which is isomorphic to space of functionals on the space of
holomorphic functions in neighborhoods of $S$ on $\Bbb CS$. Then
we can define the horospherical Cauchy transform of a functional
$f\in \operatorname {Hyp}(S)$ as its value on the function
$1/(1-\zeta\cdot z), \zeta\in \Xi_S$, which is holomorphic in a
neighborhood of $S$. It is easy to reformulate this definition in
the language of cohomology.

The horospherical Cauchy transform of cohomology $H^{(n-1)}(\Bbb
CS\backslash S, \Cal O)$ is connected with the holomorphic
cohomological language for this cohomology \cite {EGW95}. Let us
define $$ \Cal E= \{(z,\zeta); z\in \Bbb CS\backslash S, \zeta \in
\Xi_S, z\in E(\zeta)\}.$$ It is the Stein manifold. The fibering
$\Cal E\rightarrow (\Bbb CS\backslash S)$ (with the fibers - the
intersections $L(z)\bigcap \Xi_S$) satisfies the conditions in
\cite {EGW95} and we can compute the analytic cohomology of $\Bbb
CS\backslash S$ using the complex of holomorphic forms
$\phi(z,\zeta;d\zeta)$ on $\Cal E$ with the differentials only
along the fibers. If we restrict these forms on a section of the
fibering and take the $(0,q)$-part we obtain the operator on
Dolbeault cohomology. We are interested in cohomology in the
dimension $n-1$ equal to the dimension of fibers. So the closeness
of forms $\phi$ is trivial. Let us consider the operator
$$ \Cal O(\Xi_S)\rightarrow H^{(n-1)}((\Bbb CS\backslash S), \Cal
O): F(\zeta) \mapsto \phi_F(z,\zeta;d\zeta)=F(\zeta)
\nu_z(d\zeta).$$ The basic result is that this operator is an
isomorphism. More exact, if $f\in H^{(n-1)}(\Bbb CS\backslash S,
\Cal O)$ we take $F=\Cal L \hat f$. We can interpret this
construction as a holomorphic Hodge theorem: we pick up canonical
representatives in cohomology class but instead of Riemannian
geometry we use the complex one. We already used such
constructions in the context of the Penrose transform \cite
{GH78}.

In conclusion let us remark that the boundary of $\Cal E$ is the
fibering over $S$ with cycles $L_\Bbb R (x)$ as fibers. So if the
holomorphic form $\phi (z,\zeta;d\zeta)$ has boundary values in
some sense, then we can integrate the boundary form on the fibers
and receive a function on $S$. The inversion formula supports this
correspondence between regular functions and hyperfunctions.

\subhead Hyperbolic restriction \endsubhead We found that the
holomorphic inversion formula can be restricted on the compact
real form of $\Bbb CS$ - the real sphere $S$ - and gives a
possibility to reconstruct the harmonic analysis on it. It turns
out that there is a possibility to restrict this formula also on
noncompact forms. We can observe here a principal difference with
the theory of representations where we can not find
representations of real forms in such a way.

Let us start from the hyperbolic space $H^n=H$ which we will
realize as the one sheet of two-sheeted hyperboloid
$$ \square (x) = x_1^2-x_2^2-\cdots -x_{n+1}^2 =1, \qquad x_1>0.$$
It is convenient to replace the coordinates and to use here the
bilinear form corresponding to this quadratic form:
$$\zeta\cdot z=\zeta _1\cdot z_1-\zeta _2\cdot z_2\cdots -\zeta_{n+1}z_{n+1}.$$
Since $H$ is noncompact we need to put some decreasing conditions
on the class of functions. For simplicity, let us take the space
$D(H)$ of finite continuous functions. The set of complex
horospheres $E(\zeta)$ which do not intersect $H$ is parameterized
by the points of $CR$-submanifold $$\Xi_H=\{\zeta; \zeta=\lambda
\xi, \lambda \in \Bbb C\backslash \Bbb R, \xi\in \Bbb R^{n+1},
\square (\xi)=0\}.$$ We define the horospherical Cauchy transform
$\hat f(\zeta), \zeta\in \Xi_H, f\in D(H),$ by the same formula as
earlier but we replace the integration on $S$ by the integration
on $H$. The boundary is
$$ \{\partial \Xi_H=\xi\in \Bbb R^{n+1}, \square (\xi)=0\}.$$ Let
us take the boundary values $\hat f (\xi), \xi\in \partial\Xi_H,
\lambda\rightarrow 1+i0.$ In the definition of the dual
horospherical transform $\check F$ for $CR$-functions $F$ we take
the cycle $\{ \xi=x+ \mu, \mu\cdot x=0, \mu\cdot \mu=-1$\} as the
cycle $L_\Bbb R(x), x\in H,$ in the intersection-cone $L(x)\bigcap
\partial \Xi_H$. So if $x=(1,0,\dots,0)$, then $\mu_1=0,
\mu_2^2+\cdots+\mu_{n+1}^2=1$. Of course, we can replace this
cycle by any other cycle which intersects once all generators of
the cone. As it was shown in \cite {Gi00$'$} then the inversion
formula of this paper holds and it gives the possibility to
reproduce the harmonic analysis in the hyperbolic space $H$.

The difference starts when we build the analogue of the spherical
Fourier transform, The submanifold $\Xi_H$ is invariant relative
to the action noncompact Abelian subgroup $\Bbb R_+\subset \Bbb
C^\times$. The decomposition in the Mellin integral gives the
spherical irreducible representations of $SO(1,n)$ and the
composition with the horospherical Cauchy transform gives the
spherical Fourier transform. So on $\Xi_S$ the Abelian group was
compact and the spectrum was discrete; on $\Xi_H$ the Abelian id
is noncompact and the spectrum is continuous. \subhead
(2,n-1)-hyperbolic restrictions
\endsubhead We investigated the restrictions of the horospherical
Cauchy transform on the symmetric Stein space $\Bbb CS$ for 2
 Riemannian symmetric spaces $S,H$, compact and noncompact correspondingly. In the conclusion we will
discuss the restriction on one pseudo Riemannian form - the
hyperboloid $X$ of the signature (2,n-1):$$\square (x) =
x_1^2+x_2^2-x_3^2-\cdots -x_{n+1}^2 =1$$ and we will replace the
basic bilinear form by the form corresponding to this quadric
form:$$\zeta\cdot z=\zeta _1\cdot z_1+\zeta _2\cdot z_2-\zeta
_3\cdot z_3\cdots -\zeta_{n+1}z_{n+1}.$$ For $n=3$ we have
$X=SL(2;\Bbb R)$. The corresponding results on the horospherical
Cauchy transform were obtained in \cite {Gi00,Gi02}. We will see
that the picture in this case is a very interesting combination of
the pictures for $S$ and $H$.

The parametric set of complex horospheres $E(\zeta)$ which do not
intersect $X$ has 3 components: the component $\Xi_X^0$ is exactly
the same as the set $\Xi_H$ in the last example (but for another
bilinear form); 2 other the components are the connected
components of the set
$$\Xi_X^{\pm}=\{\zeta-\xi+i\eta:\square (\xi)=\square(\eta)>1\}.$$
Correspondingly, we can for $f\in D(X)$ separate 3 components
$\hat f_0(\zeta), \hat f_\pm(\zeta)$ which have these 3 sets as
the domains. They will be correspondingly $CR$-functions and
holomorphic functions. We take their boundary values. Let $F$ be a
function with such 3 components. We define the dual horospherical
transform $\check F$ of the same structure as above, but it will
have 3 components. We need to describe the sets $L_\Bbb
R^0(x),L_\Bbb R^\pm(x)$ along which we integrate in the dual
Cauchy transform. In all cases they lie in intersection of the
boundaries $\partial \Xi_X^0,\partial \Xi_X^\pm$ with $ L(x), x\in
X$. The explicit descriptions are exactly as above: $L_\Bbb
R^0(x)$ is the same as for $H$ and $L_\Bbb R^\pm(x)$ is the same
as for $S$. The difference is only in the real quadratic form. So
the hyperboloid $\{\eta: \eta\cdot x=0, \square (\eta)=1\}$ has 2
sheets. The important point is that all these 3 sets for $X$ are
not compact and we integrate along manifolds which are not cycles.
Nevertheless, they all can be compactified by the  parameters of
degenerate horospheres $$ \xi\cdot z=0, \xi \in \Bbb R^{n+1},
\square (\xi)=0, \xi\cdot x=0.$$ In such a way we construct the
cycle which has 3 components which intersect on this set. We can
define the boundary values of $\hat f$ on this set and we the
inversion formula can be written down in this case. It was proved
in \cite {Gi00}. Let us emphasize that this proof as well as
proofs in the previous cases, uses an universal inversion formula
for hyperplane sections of a quadric and the inversion formula for
horospheres is homotopic to the Radon inversion formula.

\subhead Final remarks\endsubhead 1)There is one more version of
the horospherical transform on the symmetric space $\Bbb CS^n$
which is connected with the Plancherel formula for $L^2(\Bbb
CS^n)$. So we consider consider finite functions $f\in D(\Bbb CS)$
(instead of holomorphic functions) and define the (real)
horospherical transform $$\hat f(\zeta)=\int_{E(\zeta)} f(z)
\omega \wedge \bar \omega, \qquad \zeta \in \Bbb C\Xi.$$ The dual
horospherical transform is defined as
$$\check F(z)= \int _{L(z)} F(\zeta)\nu_z(d\zeta)\wedge \overline
{\nu_z(d\zeta)}, \qquad z\in \Bbb CS.$$ There is an inversion
formula $$f=(\Cal L\bar \Cal L \hat f)^{\vee}$$ which can be
extended on $L^2(\Bbb CS)$. We consider the Fourier transform
relative to the action of $\Bbb C^\times$ on $\Bbb C\Xi$ and its
composition with the horospherical transform gives the spherical
transform (the decomposition of the representation of $SO(n+1;\Bbb
C)$ on irreducible ones). Here it will be only one series
dependent on one continuous and one discrete parameter: the
polynomial symbol of the operator $\Cal L \bar \Cal L$ is the
Plancherel density.

2) This holomorphic horospherical duality can be generalized
following the results of \cite {Gi05} on arbitrary symmetric Stein
manifolds $G/H$ where the semisimple Lie group $G$ and its
involutive subgroup $H$ are complex .

3) In the construction of the holomorphic horospherical duality we
work only with two dual families of submanifolds and do not appeal
to the groups (only the Abelian group appears if we want to have a
Fourier type transform). It would be very interesting to find some
geometrical conditions on dual families which admit similar
explicit inversion formula.


\Refs

\widestnumber \key {KKMOOT}


\ref\key {EFG95} \by M.Eastwood,S.Gindikin, H.-W.Wong \paper
Holomorphic realization of $\overline\partial$-cohomology and
constructions of representations\jour J.Geometry and Physics\vol
17\yr 1995\pages 231--244 \endref

 \ref \key {Gi00} \by S.Gindikin \paper Integral Geometry on
$SL(2;\Bbb R)$ \jour Math.Res.Letters \yr 2000 \pages 417-432
\endref

\ref \key {Gi00$'$} \by S.Gindikin\paper Integral geometry on
hyperbolic spaces\inbook Harmonic Analysis and Integral
Geometry,\ed M.Picardello\yr 2000\pages 41--46\publ Chapman and
Hall\endref


\ref \key {Gi02} \by S. Gindikin \paper An analytic separation of
series of representations for $SL(2; \Bbb R)$ \jour Moscow Math.
J. \vol 2 \issue 4 \yr 2002 \pages 1-11 \endref

\ref \key {Gi04} \by S.Gindikin \paper Complex horospherical
transform on real sphere \toappear\yr 2005\endref

\ref \key {Gi05}\by S.Gindikin \paper Horospherical Cauchy-Radon
transform on compact symmetric spaces \paperinfo preprint \yr
2004\endref

\ref \key {GH78} \by S.Gindikin and G.Henkin \paper Integral
geometry for ${\overline\partial}$-cohomology in $q$-linearly
concave domains in  $\Bbb CP^n$ \jour Funct.Anal.Appl. \lang
Russian\vol 12\issue 4 \yr 1978 \pages 6--23 \endref


 \ref \key {GH90}\by S.Gindikin and G.Henkin
\paper The Cauchy-Fantappie on projective space\jour
Amer.Math.Soc.Transl.(2)\vol 146 \yr 1990 \pages 23--32\endref

\endRefs
\enddocument